\documentclass[]{article}

\usepackage{amsthm}
\usepackage{amsmath}
\usepackage{amssymb}
\usepackage{tikz}
\usepackage{graphicx}
\usepackage[footnotesize,hang]{caption}
\usepackage{subfig}
\usepackage{color}
\usepackage{rotating}
\usepackage{mathtools}
\usepackage{xcolor}

\usepackage[margin=1.2in]{geometry}

\newcommand{\var}{\mathrm{Var}}

\usetikzlibrary{decorations.pathmorphing}

\renewcommand{\inf}{\infty}

\newcommand{\e}{\mathrm{e}}

\renewcommand{\d}{\,\mathrm{d}}
\newcommand{\diff}[2]{\frac{\mathrm{d} #1}{\mathrm{d} #2}}
\renewcommand{\max}{\mathrm{max}}
\renewcommand{\min}{\mathrm{min}}

\def\XXint#1#2#3{{\setbox0=\hbox{$#1{#2#3}{\int}$}
\vcenter{\hbox{$#2#3$}}\kern-.5\wd0}}

\newtheorem{theorem}{Theorem}[section]
\newtheorem{example}{Example}[section]

\theoremstyle{remark}
\newtheorem*{remark}{Remark}

\title{On Asymptotics of Optimal Stopping Times}
\author{C. J. Lustri\footnote{Electronic address: christopher.lustri@mq.edu.au}, \ G. Yu. Sofronov\footnote{Electronic address: georgy.sofronov@mq.edu.au} \ and H. N. Entwistle\footnote{Electronic address: hugh.entwistle@mq.edu.au}}
\date{%
Department of Mathematics and Statistics, 12 Wally's Walk, Macquarie University, New South Wales 2109, Australia
}                                     

\begin{document}
\maketitle

\begin{abstract}
We consider optimal stopping problems, in which a sequence of independent random variables is drawn from a known continuous density. The objective of such problems is to find a procedure which maximizes the expected reward. 
In this analysis, we obtain asymptotic expressions for the expectation and variance of the optimal stopping time as the number of drawn variables becomes large. In the case of distributions with infinite upper bound, the asymptotic behaviour of these statistics depends solely on the algebraic power of the probability distribution decay rate in the upper limit. In the case of densities with finite upper bound, the asymptotic behaviour of these statistics depends on the algebraic form of the distribution near the finite upper bound. Explicit calculations are provided for several common probability density functions, which are compared to numerical simulations that support the asymptotic predictions.
\end{abstract}

%
%



\section{Introduction}

Optimal stopping problems pose the challenge of deciding when to stop some stochastic process in order to maximise some objective, or utility. This problem has arisen in a number of contexts, such as deciding how many candidates for a job should be interviewed before deciding upon one particular candidate, known in literature as the secretary problem \cite{Ferg89,SEALE1997221,STEIN2003140}. Optimal stopping has been used to determine when assets should be bought and sold in order to maximize profits; this has been described as the house-selling problem, and is discussed in \cite{DAVID1998576,karlin1962stochastic,multiplesellSofronov}. Optimal stopping has been incorporated into the theory of online auctions, as in \cite{hajiaghayi2004adaptive}. In this study, the authors link the analysis of the secretary problem to the design of a mechanism for online auctions. Related analyses may be found in \cite{harrell2015online,kleinberg2005multiple}. Optimal stopping also features in financial and economic applications such as in the pricing of American options \cite{egloff2005monte}, games \cite{Ivashko2020}, operational risk insurance \cite{Targino2017}, dynamic pricing \cite{karpowicz2007double}, and more complicated buying-selling problems \cite{Sofronov2006,Sofronov2016,Sofronov_MCAP2020}.



Optimal stopping problems are formulated in terms of observing random variables, and determining the stopping point in order to maximize a particular reward function. The problem considered here involves observing a sequence of random variables $y_1,y_2,\dots,y_N$, and making the decision to stop after a particular number of observations, denoted $m$ where $1 \leq m \leq N$, based on the variables that have been previously observed at that stage. After stopping, we receive a reward which is a function of the already observed values $y_1,\ldots,y_m$. 
This problem is one representative of  a class of optimal stopping problems that consists of finding a sequential procedure that maximizes the expected reward. For a more extensive discussion of this class of problem, see \cite[Section 13.4]{degroot2005optimal}). 

While there exists extensive literature on the theory of optimal stopping problems \cite{Chow,degroot2005optimal}, less attention has been paid to asymptotic properties of stopping times. Most existing asymptotic results focus on ``no-information'' problems in which the distribution of the observations is unknown.  


In no-information problems, an observer determines the relative rank of each observation, and the reward function is a function of these ranks. An example of a no-information problem is the secretary problem \cite{Ferg89,freeman1983}. In this problem, the reward function is the indicator variable of the best object, which means that the observer aims to maximise the probability that the best object is selected. As a consequence, the secretary problem is sometimes referred to as the no-information best-choice problem. 

There are many other variations of the no-information problem and significant work exists on the asymptotic properties of the stopping time (see, for example, \cite{demers2018note,freeman1983,UnifiedApproach}). It was shown in \cite{mazalov2004asymptotic} that in the secretary problem with a sequence of $N$ observations, the asymptotic expectation and variance for the stopping time is ${2N}/{\e}$ and $\left( {2}/{\e} - {5}/{\e^2} \right)N^2  \approx 0.059 N^2$, respectively. Asymptotic descriptions of statistical properties for other no-information problem variants can be found in \cite{demers2018note,GilbMost},  
where the techniques used are dependent on the structure and variation of each problem. The asymptotics where the number of observations is random has also been studied in  \cite{bestchoice_randomN}.

There is substantially less literature describing the asymptotic behaviour of statistical properties for ``full-information" problems,  when the distribution of random variables is known beforehand. Gilbert and Mosteller \cite{GilbMost} studied the optimal stopping strategy for the full-information problem in which the objective is identical to that of the secretary problem, known as the full-information best-choice problem \cite{Gnedin1996fullinfo}. In this setting, the optimal stopping strategy was shown to consist of stopping and choosing the $m$-th observation $y_m$ if it is the highest ranked out of all observations made at this point, and has a value exceeding a threshold, which depends on $m$. 
Asymptotically, it was shown that $v_n$, the value of a sequence with $n$ steps remaining, $v_n \sim 1-c/n$, where $c\approx 0.804352$. 

In best-choice settings, variables that do not follow the uniform distribution can be re-scaled by applying the cumulative density function of the particular distribution to achieve a uniform distribution that is monotonically equivalent as far as the best-choice is concerned. The asymptotic probability of finding the best object (see  \cite[Section 3]{GilbMost} and \cite{samuelsfullinfo}) is independent of the distribution of the variables, tending to approximately 0.58 as the number of observations in the sequence becomes large.

The best-choice problems that have been described so far have been interested in finding the very best object. For the full-information case, this condition can be relaxed by instead setting our reward at stopping time $m$ to $y_m$ and seeking to maximise the expected reward. A special case of this problem called the uniform game \cite[Section 5a]{GilbMost} is related to the well-known Cayley's problem (see, for example, \cite{Ferg89,Moser}). In \cite{GilbMost}, the authors showed that the expected reward of a sequence of $n$ independent and identically distributed (iid) random variables having the standard uniform distribution can be approximated in the following way (see also \cite{Moser})
\begin{equation}
v_n\sim 1-\frac{2}{n+\log(n+1)+1.767} \qquad \mathrm{as} \qquad n \rightarrow \infty.\label{v_n_GilbMost} 
\end{equation}
In \cite{mazalov2004asymptotic}, the authors found the asymptotics of the expected value and the variance of the stopping time as $N/3$ and $N^2/18$ respectively for the uniform game. The asymptotic techniques used to compute this behaviour are highly dependent on the distribution, and cannot be easily generalized to more general classes of distribution. In \cite{LimitTheorems1990,kennedy1991asymptotic}, using the extreme value theory, the authors prove limit theorems for threshold-stopped random variables and derived the asymptotic distribution of the reward sequence of the optimal stopping (iid) random variables.  


In the full-information problem where the reward at time $m$ is $y_m$, such as that considered in the present study, the distribution of the observations plays a significant role in the asymptotic behaviour of the outcome statistics. This study offers a general asymptotic technique for calculating the asymptotic behaviour of $E(\tau_N)$ and $\var(\tau_N)$ as $N\to\infty$ for general classes of probability distributions in the full-information problem where we wish to maximise the expected reward $y_m$. In addition to recovering existing results from \cite{mazalov2004asymptotic} for the uniform distribution, we extend the analysis to consider the effect of drawing elements of the sequence from a wide range of common distributions.


\section{Formulation}
Let $y_1,y_2,\dots,y_N$ be a sequence of independent random variables with $y_m$ as the reward at time $m$.
This problem of finding the optimal stopping rule can be solved by backward induction using the following recurrent equation (see, for example, \cite{Chow,degroot2005optimal}):
\begin{equation}\label{1:vn}
v_n = E(\max\{y_{N-n+1},v_{n-1}\}),\quad n=1,\dots,N,\; v_0=-\infty,
\end{equation}
where $v_n$ is the value of a sequence with $n$ steps, $v_N$ is the expected reward. Here $\max\{y_{N-n+1},v_{n-1}\}$
represents the maximum gain that is possible to obtain having $n$ steps left. If $n = 1$,
we have to stop and our gain will be $y_N$. If $1< n \leqslant N$,
we can either stop or continue. If we stop, our gain is $y_{N-n+1}$, and if we continue, our expected gain
is $v_{n-1}$.

If $y_1,y_2,\dots,y_N$ are iid continuous random variables with common probability density function (pdf) $f(y)$, then
\begin{equation}\label{value}
 v_n=\int_{-\infty}^{\infty} \max\{y,v_{n-1}\} f(y)\, dy,\quad n=1,\dots,N,\; v_0=-\infty.
\end{equation}

The optimal stopping rule $\tau_N$ is
\begin{equation}\label{1:opstop}
\tau_N =\min \{ m : 1\leqslant  m \leqslant N, y_m\geqslant v_{N-m}\}.
\end{equation}

This optimal stopping rule formulation is for single stopping. It can be extended for multiple stopping for sums of random variables, see \cite{SofronovMultStop,Sofronov2013}. 

\section{Computing $v_n$ behaviour}\label{S:Inf}

\subsection{Recurrence relation for $v_n$}

Assume that $y_n$ is drawn from a continuous pdf $f(y)$, which has unbounded support in the positive direction. This function has a cumulative distribution function (cdf) $F(y)$, and a  complementary cumulative distribution function (ccdf) $h(y) = 1-F(y)$. 




\begin{theorem}\label{1 Thm Inf}
Let $Y$ be a random variable whose expectation exists, and which is drawn from a continuous distribution $f(y)$. Assume that there exists positive $\lambda$ and $\Delta$ such that the ccdf satisfies $h(y) \leq \lambda/y^{1+\Delta}$ for sufficiently large $y$. 

The value of a sequence with $n+1$ steps is given by
\begin{equation}
v_{n+1}= v_n +  \int_{v_n}^{\infty} h(y) \mathrm{d}y. \label{2:vn recur inf}
\end{equation}

\end{theorem}
\begin{proof}
The recurrence relation \eqref{1:vn} can be written as 
\begin{align}
v_{n+1} &= \int_{-\infty}^{v_n} v_n f(y) \mathrm{d}y + \int_{v_n}^{\infty} y f(y) \mathrm{d}y \\&= v_n F(v_n) + \lim_{K\rightarrow \infty} \left\{[y F(y)]_{v_n}^K - \int_{v_n}^K F(y) \mathrm{d} y\right\} =  v_n + \lim_{K\rightarrow \infty} \left\{K h(K) + \int_{v_n}^K h(y) \mathrm{d} y\right\}.
\end{align}
From our assumptions, it is easy to see that $\lim_{K\rightarrow\infty} Kh(K) = 0$, and that the integral must converge in the limit that $K \rightarrow \infty$. This therefore gives \eqref{2:vn recur inf}, and completes the proof.
\end{proof}
%

We note that if $f(y)$ has bounded support in the positive direction such that $f(y) = 0$ for $y > y_{\mathrm{max}}$, a similar recurrence relation may be obtained in nearly identical fashion.
\begin{theorem}
Let $Y$ be a random variable whose expectation exists, and which is drawn from a continuous distribution $f(y)$ with bounded support in the positive direction, such that $f(y) = 0$ for $y > y_{\mathrm{max}}$. Assume that there exists positive $\lambda$ and $\Delta$ such that the ccdf satisfies $h(y) \leq \lambda/y^{1+\Delta}$ for sufficiently large $y$. 

The value of a sequence with $n+1$ steps is given by
\begin{equation}
v_{n+1}= v_n +  \int_{v_n}^{y_{\mathrm{max}}} h(y) \mathrm{d}y. \label{2:vn recur fin}
\end{equation}
\end{theorem}
\begin{proof}
For $y > y_{\mathrm{max}}$, it is clear from the definition of the ccdf that $h(y) = 0$. Therefore
\begin{equation}
\int_{v_n}^{\infty} h(y) \d y = \int_{v_n}^{y_{\mathrm{max}}} h(y) \d y +  \int_{y_{\mathrm{max}}}^{\infty} h(y) \d y = \int_{v_n}^{y_{\mathrm{max}}} h(y) \d y.
\end{equation}
Replacing the integral appropriately in Theorem \ref{1 Thm Inf} completes the proof.
\end{proof}

\subsection{Asymptotics of $v_n$ as $n \rightarrow \infty$}

Using the asymptotic controlling factor method (found in Section 5.3 of \cite{bender2013}), we know that if a function $v_n$ grows no more rapidly than $\e^{a n^b}$ with $b < 1$ as $n \rightarrow \infty$, then the leading-order asymptotic solution for $v_n$ satisfies
\begin{equation}
v_{n+1} - v_n \sim v_n' \qquad \mathrm{as} \qquad n \rightarrow \infty,
\end{equation}
where asymptotic equivalence is defined in the usual fashion, and $'$ denotes differentiation with respect to $n$. It is straightforward to obtain a finite upper bound for the integral in \eqref{2:vn recur inf}, which implies that the growth of $v_n$ can be no faster than linear in $n$. Hence, this theorem holds. 

We can therefore use \eqref{2:vn recur inf} to write the asymptotic relation 
\begin{equation}\label{2:vn ode}
v_n' \sim  \int_{v_n}^{\infty} h(y) \mathrm{d}y  \qquad \mathrm{as} \qquad n \rightarrow \infty,
\end{equation}
with the upper bound replaced by $b$ in the finite upper support case, corresponding to \eqref{2:vn recur fin}. 

\begin{remark}
The expression in \eqref{2:vn ode} can be manipulated to obtain a convenient asymptotic representation for $v_n$, subject to some additional assumptions. Suppose that $v_n''$ exists, is integrable for sufficiently large $n$, and has an asymptotic series representation. These conditions are sufficient to allow for both sides of \eqref{2:vn ode} may be differentiated (see, for example, the discussion in Section 3.8 of \cite{bender2013}) to give
\begin{equation}\label{2:vddn}
v''_n \sim -h(v_n)v'_n  \qquad \mathrm{as} \qquad n \rightarrow \infty.
\end{equation}
Assume that there exists some $k$ such that $v_n' \neq 0$ for $n > k$. Equation~\ref{2:vddn} therefore gives   
\begin{equation}\label{2:hvn}
h(v_n) \sim - \frac{v''_n}{v'_n} \qquad \mathrm{as} \qquad n \rightarrow \infty.
\end{equation}
The assumptions required for \eqref{2:hvn} to be valid hold for all distributions considered in the present study. This relationship will therefore be used to simplify a later result for a collection of common distributions.
\end{remark}

%
%
\subsection{Example Calculations}
We illustrate the application of these ideas to some common distributions, such as the exponential distribution, which is studied in Example \ref{E1}. In Example \ref{E2}, we demonstrate that the same method can be applied to heavy-tailed distributions by considering the Pareto distribution, and in Example \ref{E3} we show that this method can easily be applied to distributions with finite support by considering the uniform distribution. 

Finally, in Example \ref{E4} we consider a broad class of distributions with exponential upper tails and show that the asymptotic behaviour of $v_n$ is fully determined for this class of distribution, subject to the assumption that $v_n$ increases without bound.

\begin{example}\label{E1}
The exponential distribution is given by
\begin{equation}
f(y) = \frac{1}{\beta}\e^{-y/\beta},
\end{equation}
with $y \in [0,\inf)$ and $\beta > 0$.
\end{example}

The ccdf $h(y)$ is given by $h(y) = \e^{-y/\beta}$. Hence, \eqref{2:vn ode} becomes
\begin{equation}
\diff{v_n}{n} \sim \beta \e^{-v_n/\beta}  \qquad \mathrm{as} \qquad n \rightarrow \infty.
\end{equation}
Solving this formal equation gives the large-$n$ asymptotic behaviour as
\begin{equation}
v_n \sim \beta \log(n) \qquad \mathrm{as} \qquad n \rightarrow \infty.
\end{equation}
Note that the constant term obtained by solving the ordinary differential equation does not contribute to the leading-order expression in the large-$n$ limit.

\begin{example}\label{E2}
The Pareto distribution is given by
\begin{equation}
f(y) = \frac{\alpha\beta^\alpha}{y^{\alpha+1}},
\end{equation}
on $y \in [\beta,\infty)$, where $\beta > 0$. We assume that $\alpha > 1$.
\end{example}
The ccdf is given by $h(y) = (\beta/y)^{\alpha}$. The differential equation for the large-$n$ asymptotic behaviour is therefore given by
\begin{equation}
\diff{v_n}{n} \sim \left(\frac{\beta^{\alpha}}{\alpha-1}\right)v_n^{1-\alpha}.
\end{equation}
Solving this formal equation gives the asymptotic behaviour in the large-$n$ limit as
\begin{equation}
v_{n} \sim  \beta\left(\frac{\alpha n}{\alpha-1}\right)^{1/\alpha}  \quad \mathrm{as} \quad n \rightarrow \infty.
\end{equation}

\begin{example}\label{E3}
The uniform distribution is given by
\begin{equation}
f(y) =\frac{1}{b-a},
\end{equation}
on $y \in [a,b]$, where $b > a$.
\end{example}
The ccdf is given by $h(y) =\tfrac{b-y}{b-a}$.  The differential equation for the large-$n$ asymptotic behaviour is therefore given from \eqref{2:vn recur fin} by
\begin{equation}
\diff{v_n}{n} \sim  \frac{(b-v)^2}{2(b-a)} \qquad \mathrm{as} \qquad n \rightarrow \infty.
\end{equation}
Solving this formal equation gives the asymptotic behaviour in the large-$n$ limit as
\begin{equation}
v_{n} \sim  b - \frac{2(b-a)}{n}  \quad \mathrm{as} \quad n \rightarrow \infty.
\end{equation}

\begin{example}\label{E4}
A distribution $f(y)$ is given with a ccdf that, for sufficiently large $y$, satisfies 
\begin{equation}\label{2.bound}
\left| h(y) - \gamma  \e^{-(y/\beta)^{\alpha}}\right| < \frac{\e^{-(y/\beta)^{\alpha}}}{y^{\Delta}},
\end{equation}
for positive $\Delta$, and where $\alpha$, $\beta$, and $\gamma$ are positive constants. Assume that $v_n$ increases without bound (ie. for any $N$, there exists an $n$ such that $v_n > N$).
\end{example}
The condition in \eqref{2.bound} gives the asymptotic behaviour of $h(y)$ for large $y$. Assuming $v_n$ increases without bound as $n$ increases, there must be some sufficiently large $n$ such that this inequality holds. We may therefore replace the full expression for $h(y)$ in \eqref{2:vn ode} with its asymptotic value in the large-$n$ limit, giving the ordinary differential equation
\begin{align}
\diff{v_n}{n}  \sim  \int_{v_n}^{\infty} h(y)\d y = \frac{\gamma \beta}{\alpha}\Gamma\left(\frac{1}{\alpha},\frac{v_n^{\alpha}}{\beta^{\alpha}}\right) + g(n),\label{2.ode4}
\end{align}
where $\Gamma$ represents the upper incomplete gamma function, and 
\begin{equation}
|g(n)| < \int_{v_n}^{\infty}  \frac{\e^{-(y/\beta)^{\alpha}}}{y^{\Delta}} \d y < \frac{1}{v_n^{\Delta}}\int_{v_n}^{\infty}\e^{-(y/\beta)^{\alpha}}\d y.
\end{equation}
As $v_n$ is assumed to increase without bound, this integral is asymptotically subdominant compared to the right-hand side of \eqref{2.ode4} as $n \rightarrow \infty$. The solution to the differential equation \eqref{2.ode4} may therefore be approximated using standard asymptotic methods to give the formal asymptotic relation
\begin{equation}
v_{n} \sim \beta\log(n)^{1/\alpha} \quad \mathrm{as} \quad n \rightarrow \infty.
\end{equation}
It is straightforward to show using \eqref{2:hvn} that large-$n$ asymptotics of $h(v_n)$ are independent of $\alpha$ and $\beta$ to leading order, giving $h(v_n) \sim 1/n$ in the limit that $n \rightarrow \infty$.

\subsubsection{Tabulation of Further Examples}

Using similar methods to the previous examples, the asymptotic behaviour of $v_n$ may be computed for a wide range of common distributions. Table~\ref{T:vn} contains the large-$n$ asymptotics for a number of common distributions.




These densities contain a number of parameters, some of which are required to satisfy particular conditions. The normal distribution permits arbitrary $\mu$, but requires $\sigma > 0$. The gamma distribution requires $\alpha > 0$ and $\beta > 0$, and the function $\gamma$ represents the lower incomplete gamma function, found in \cite{NIST:DLMF}. We also assume that $\alpha \neq 1$ for the purposes of the calculation shown here. The $\alpha = 1$ case requires a different analysis, and corresponds to the exponential distribution. The Pareto distribution requires $\beta > 0$ and $\alpha > 1$.

The uniform distribution requires $a < b$, while the triangular  distribution requires $a < c < b$. The Wigner distribution requires $R > 0$, and the beta distribution requires $\alpha > 0$ and $\beta > 0$. In the expression for the beta distribution, the function $B(\alpha, \beta)$ represents the standard beta function, while $B(y; \alpha, \beta)$ represents the incomplete beta function (see \cite{NIST:DLMF}).

The function $W$ denotes the Lambert-W function. The asymptotic expressions for $v_n$ associated with the normal and gamma distributions can be further simplified by observing that the argument of the Lambert-W function becomes large in the asymptotic limit, and
\begin{equation}\label{lambert}
W(x) = \log(x) - \log(\log(x)) + \frac{\log(\log(x))}{\log(x)} + {o}\left(\frac{1}{\log(x)}\right)\quad \mathrm{as} \quad x \rightarrow \infty.
\end{equation}
Using the first two terms of this expansion produces an asymptotic expression for $v_n$ with error that is $o(1)$ as $n \rightarrow \infty$ for both the normal and gamma distributions.

We note that the results obtained using this asymptotic formulation for the uniform distribution are consistent with previous analyses from \cite{demers2018note, mazalov2004asymptotic}.

Importantly we see that $v_n$ approaches a maximum value for each distribution with compact support, corresponding to the maximum possible value of $y$ in the domain, as assumed and subsequently confirmed in the previous analysis. In contrast, for each pdf on a domain with unbounded upper support, $v_n$ increases monotonically without bound. The asymptotic behaviour of $v_n$ will be subsequently used to determine the expectation and variance of $\tau_N$ for each example.


\begin{table}[ht]
  \begin{center}
    \caption{This table details the behaviour of $v_n$ in the large-$n$ limit for several common probability density functions. For each example, the table contains the pdf equation $f(y)$ and the density domain. The next column contains the ccdf $h(y)$. The final column shows the asymptotic behaviour of $v_n$ in the limit that $n \rightarrow \infty.$ In each case, $v_n$ increases without bound.}\label{T:vn}

    \label{Tab:Recurrence}
    \begin{tabular}{|l|l|l|} 
        \hline
      Distribution   & $h(y)$ & Tail asymptotics as $n \rightarrow \infty$\\
      \hline
                  Normal; $y \in (-\infty,\infty)$: &&\\[0pt]
                $f(y) = \frac{1}{\sqrt{2\pi \sigma^2}} \mathrm{e}^{-(y-\mu)^2/2\sigma^2}$  & $\frac{1}{2} - \frac{1}{2}\mathrm{erf}\left(\frac{y-\mu}{\sigma \sqrt{2}}\right)$ & $ v_n \sim \mu + \sigma \sqrt{W\!\left(\frac{n^2}{2\pi}\right)}$\\[10pt]\hline
                            Gamma; $y \in (0,\infty)$: &&\\[0pt]
                                                              $f(y) = \frac{\beta^{-\alpha}}{\Gamma(\alpha)}y^{\alpha-1}\mathrm{e}^{- y/\beta}$  & $1-\frac{\gamma(\alpha, y/\beta)}{\Gamma(\alpha)}$ & {$\!\begin{aligned}  v_n \sim \beta&(\alpha-1) \,W\!\left[\tfrac{1}{\alpha - 1}\left(\tfrac{n}{\Gamma(\alpha)}\right)^{1/(\alpha-1)}\right]\end{aligned}$}\\[10pt]\hline
                       Triangular; $y \in [a,b]$: && \\[5pt]
                        {$\!\begin{aligned} y \leq c \,:& \tfrac{2(y-a)}{(b-a)(c-a)}\\  
                                                      y > c \,:& \tfrac{2(b-x)}{(b-a)(b-c)}\end{aligned}$} &{$\!\begin{aligned} y \leq c \,:& \tfrac{(y-a)^2}{(b-a)(c-a)}\\  
                                                      y > c \,:& \tfrac{(b-x)^2}{(b-a)(b-c)}\end{aligned}$} & $v_n \sim b-\sqrt{\frac{3(b-a)(b-c)}{2n}}$ \\[15pt]\hline
                         Wigner; $y \in [-R,R]$: &&  \\[5pt]
  $f(y) = \frac{2\sqrt{R^2 - y^2}}{\pi R^2}$ &
               $\tfrac{1}{2} - \tfrac{x\sqrt{R^2 - y^2}}{\pi R^2} - \tfrac{1}{\pi}\mathrm{arcsin}\left(\tfrac{y}{R}\right) $ & $ v_n \sim R - \frac{1}{2}\left(\frac{5\pi}{2n}\right)^{2/3}$ \\[10pt]\hline
            Beta; $y \in (0,1)$:  && \\[5pt]
                        $f(y) = \frac{y^{\alpha-1}(1-y)^{\beta-1}}{B(\alpha,\beta)}$  &$1-\frac{B(y;\alpha,\beta)}{B(\alpha,\beta)}$ & $ v_n \sim 1 - \left(\frac{B(\alpha,\beta)}{(\beta+1) n}\right)^{1/\beta}$  \\[10pt]

      \hline
    \end{tabular}
  \end{center}
\end{table}

\section{Calculating optimal stopping statistics}\label{s.opstat}

\subsection{Calculating the expectation}

Let $w_i = P(y < v_i) = 1 - h(v_i)$. The expectation is now given by
\begin{align}
E(\tau_N) &= \sum_{n=1}^N n P\\
&= \sum_{n=1}^N n P(y_1<v_{N-1},\dots,y_{n-1}<v_{N-n+1},y_n\geqslant v_{N-n})\\
&= (1-w_{N-1}) + 2 w_{N-1}(1-w_{N-2}) + \ldots + N w_{N-1}\ldots w_1\\
&=1 + \sum_{n=1}^{N-1}\prod_{i=n}^{N-1} w_i.
\end{align}
We select a value $k$ such that $0 \ll k \ll N$; a representative choice is $k = \lfloor\sqrt{N}\rfloor$. This choice does not affect the asymptotic values computed here, but would impact the form of higher order corrections.  We write the sum of products as
\begin{equation}\label{3:E}
E(\tau_N) = 1 +  \sum_{n=1}^{k-1}\prod_{i=n}^{N-1} w_i + \sum_{n=k}^{N-1}\prod_{i=n}^{N-1} w_i.
\end{equation}

We now consider the first summation term in \eqref{3:E}. Recalling that $w_n$ is non-negative and bounded above by one, we see that
\begin{equation}\label{3:ksum}
0 < \sum_{n=1}^{k-1}\prod_{i=n}^{N-1} w_i  < \sum_{n=1}^{k-1}\prod_{i=n}^{N-1} 1 = k-1.
\end{equation}


The purpose of \eqref{3:E} is to split the sum such that an asymptotic expression can be used for $w_i$ in order to approximate the second summation term. As $k$ is large in the limit that $N \rightarrow \infty$, we may use the asymptotic approximation for $v_n$ in the large $n$ limit for each term in the product. Consequently, if $v_n$ satisfies the requirements for \eqref{2:hvn} to hold, we find that
\begin{align}
\sum_{n=k}^{N-1}\prod_{i=n}^{N-1} w_i \sim \sum_{n=k}^{N-1}\prod_{i=n}^{N-1} \left(1 + \frac{v''_i}{v'_i}\right)  \quad \mathrm{as} \quad N \rightarrow \infty.\label{3:Nsum0}
\end{align}
In fact, for many distributions, this may be further simplified, using the large-$n$ asymptotics for $v_n$ and its derivatives. For the distributions in Examples \ref{E1}--\ref{E4}, as well as each distribution from Table~\ref{T:vn}, it is straightforward to show by direct calculation using \eqref{2:hvn} that 
\begin{equation}\label{2:vratio}
    h(v_n) \sim \frac{\lambda}{n} \quad \mathrm{as} \quad n \rightarrow \infty,
\end{equation}
for some positive constant $\lambda$. In fact, from Example~\ref{E4}, we see that any distribution with an exponential upper tail will have $\lambda = 1$. For distributions where $v_n$ satisfies \eqref{2:vratio},
\begin{align}
\sum_{n=k}^{N-1}\prod_{i=n}^{N-1} w_i &\sim \sum_{n=k}^{N-1}\prod_{i=n}^{N-1} \left(1 - \frac{\lambda}{i}\right) = \sum_{n=k}^{N-1} \frac{(n-\lambda)_{N-n}}{(n)_{N-n}} \sim \frac{N-\lambda-1}{\lambda+1}\sim \frac{N}{\lambda+1}.\label{3:Nsum} 
\end{align}
in the limit that $N \rightarrow \infty$, where $(a)_n$ is the Pochhammer symbol \cite{NIST:DLMF}. 
By comparing \eqref{3:Nsum} with \eqref{3:ksum}, we see that the second summation term dominates the expression for large $N$. Hence, 
\begin{equation}\label{3:Ea}
E(\tau_N) \sim  \frac{N}{\lambda + 1} \quad \mathrm{as} \quad N \rightarrow \infty.
\end{equation}
We note that, while this leading order approximation is correct in the limit that $N \rightarrow \infty$, for practical purposes it is often necessary to retain the first correction term in the approximation for $w_n$ in the limit that $n \rightarrow \infty$. While this does not give an expression as simple as \eqref{3:Ea}, it does significantly increase the accuracy of the approximation for moderately large values of $N$.

For example, in Section \ref{s.num}, we compute the expectation (and variance) of $\tau_N$ for the normal distribution. Computing only the leading-order expression for the expectation and variance of $\tau_N$ causes this computation to converge to the asymptotic value extremely slowly, due to the $\log(\log(n))$ terms present in the first correction to $w_n$ as $n \rightarrow \infty$. By including these terms in the asymptotic expansion for $w_n$, we are able to show that the simulations agree with the asymptotic predictions within the simulated range of $N$.

\subsection{Calculating the variance}

A similar process may be used to determine the square expectation, and hence the variance. The square expectation is given by
\begin{align}
E(\tau_N^2) &= \sum_{n=1}^N n^2 P\\
&= \sum_{n=1}^N n^2 P(y_1<v_{N-1},\dots,y_{n-1}<v_{N-n+1},y_n\geqslant v_{N-n})\\
&= (1-w_{N-1}) + 2^2 w_{N-1}(1-w_{N-2}) + \ldots + N^2 w_{N-1}\ldots w_1\\
&=1 + \sum_{n=1}^{N-1}(2N+1-2n)\prod_{i=n}^{N-1} w_i.
\end{align}
We again define $k$ such that $0 \ll k \ll N$, and split this series to obtain
\begin{equation}\label{4:E2}
E(\tau_N^2) = 1 +  \sum_{n=1}^{k-1}(2N+1-2n)\prod_{i=n}^{N-1} w_i + \sum_{n=k}^{N-1}(2N+1-2n)\prod_{i=n}^{N-1} w_i.
\end{equation}
We now consider the first summation term in \eqref{4:E2}. As before, we note that $w_i$ is positive and less than one, giving
\begin{equation}\label{4:FiniteVarSum}
0 < \sum_{n=1}^{k-1}(2N+1-2n)\prod_{i=n}^{N-1} w_i < \sum_{n=1}^{k-1}(2N+1-2n)\prod_{i=n}^{N-1} 1 = 2kN - 2N-k^2+2k-1.
\end{equation}

We have again split the sum such that an asymptotic expression can be used for $w_i$ in order to approximate the second summation term using the asymptotic approximation for $v_n$ in the large $n$ limit for each term in the product. If $v_n$ satisfies the requirements for \eqref{2:hvn} to hold, we find that
\begin{align}
\sum_{n=k}^{N-1}(2N+1-2n)\prod_{i=n}^{N-1} w_i \sim \sum_{n=k}^{N-1}(2N+1-2n)\prod_{i=n}^{N-1} \left(1+ \frac{v''_i}{v'_i}\right) \quad \mathrm{as} \quad N \rightarrow \infty.
\end{align}
We again observe that for many distributions, including those in Examples \ref{E1}--\ref{E4}, and each distribution from Table~\ref{T:vn}, the asymptotic behaviour of $v_n$ and its derivatives may be used to simplify $h(v_n)$, giving the expression in \eqref{2:vratio}. Noting this result, the sum can be evaluated to give
\begin{align}
\sum_{n=k}^{N-1}(2N+1-2n)\prod_{i=n}^{N-1} w_i &\sim \sum_{n=k}^{N-1}(2N+1-2n)\prod_{i=n}^{N-1} \left(1- \frac{\lambda}{i}\right)\\
&= \sum_{n=k}^{N-1} \frac{(2N+1-2n)(n-\lambda)_{N-n}}{(n)_{N-n}}\\
 &\sim \frac{(2N + \lambda+2)(N-\lambda-1)}{(\lambda+1)(\lambda+2)} \sim \frac{2N^2}{(\lambda+1)(\lambda+2)}
\end{align}
in the limit that $N \rightarrow \infty$. This expression dominates the sum in \eqref{4:FiniteVarSum}. Hence, the evaluating the asymptotic behaviour of the variance as $N \rightarrow \infty$ gives
\begin{equation}
\mathrm{Var}(\tau_N) = E(\tau_N^2) - E(\tau_N)^2 \sim \frac{2N^2}{(\lambda+1)(\lambda+2)} - \left(\frac{N}{\lambda+1}\right)^2.
\end{equation} 
This expression may be simplified to give an asymptotic approximation for the variance,
\begin{equation}
\label{4:Var}
\mathrm{Var}(\tau_N) \sim \frac{\lambda N^2}{(\lambda+1)^2(\lambda+2)} \quad \mathrm{as} \quad N \rightarrow \infty.
\end{equation}
As before, it is sometimes of practical value to retain higher corrections in the approximation for $w_n$ as $n \rightarrow \infty$, such as in the computations performed in Section [REF] on the normal distribution. This allows for the asymptotic prediction to be compared with the simulations for the moderately large values of $N$ considered here.

\subsection{Example calculations}

In Table~\ref{T:eqvar}, we compute the expectation and variance for the common pdfs computed in Examples \ref{E1}--\ref{E3}, and those shown in Table~\ref{T:eqvar}. This table illustrates the asymptotics of $h(y)$ in the limit that $ y \rightarrow \infty$ for pdfs with infinite upper support, and in the limit that $y \rightarrow y_{\mathrm{max}}$ for domains with finite upper support. For each pdf, the expectation and variance of $\tau_N$ are obtained using the expectation formula from \eqref{3:Ea} and the variance formula from \eqref{4:Var}.


In each of the densities with exponentially decaying tails (exponential, normal, and gamma), the expectation and variance are identical, corresponding to setting $\lambda = 1$ in the expectation and variance formulae. This is consistent with the result of Example \ref{E4}, which showed that the large-$n$ asymptotics are identical for any pdf with an exponentially decaying upper tail, corresponding to $\lambda = 1$. 

In the remaining pdfs, however, the expectation and variance are not identical. For distributions with heavy tails or finite upper bound, the optimal stopping statistics are determined by the rate which the distribution function decays as the upper bound is approached, characterised by the asymptotic behaviour of $h(y)$ in this limit. Consequently, the expectation and variance associated with each of these distributions take different values.

We can see that setting $\beta = 1$ in the beta distribution gives identical expectation and variance to the uniform distribution, while setting $\beta = 2$ or $\beta = 3/2$ give identical expectation and variance to the triangular and Wigner distributions respectively. This is caused by the fact that the asymptotic decay of $h(y)$ in the limit $y \rightarrow y_{\mathrm{max}}$ occurs at the same algebraic power in each case, and these statistics depend entirely on the algebraic power of the decay rate.

\begin{table}[ht]
  \begin{center}
    \caption{This table contains the asymptotic behaviour of optimal stopping statistics for several common probability density functions, as well as intermediate quantities used to compute these statistics. For each distribution, the table describes the pdf name and domain, the asymptotic behaviour of $h(y)$ in the limit that $y \rightarrow \infty$ or $y \rightarrow y_{\mathrm{max}}$ for domains with infinite and finite upper support respectively, and the value of $\lambda$ associated with this asymptotic behaviour. The final two columns contain the expectation and variance of $\tau_N$ in the limit that $N \rightarrow \infty$.}\label{T:eqvar}
    \begin{tabular}{|l|c|c|c|c|c|} 
        \hline
      Distribution & Domain & $h(y)$ as $y \rightarrow \infty$& $\lambda$ & $E(\tau_N)$ & Var($\tau_N$)\\
      \hline
            &&&&&\\
      Exponential & $[0,\infty)$ &  $ \frac{1}{\beta} \mathrm{e}^{- y/\beta}$ &$1$ & $\frac{N}{2}$ &$ \frac{N^2}{12} $ \\[10pt]
            Normal & $(-\infty,\infty)$ & $\frac{\sigma}{x\sqrt{2\pi}} \mathrm{e}^{- (y-\mu)^2/2\sigma^2}$ &$1$ & $\frac{N}{2}$ &$ \frac{N^2}{12} $ \\[10pt]
                      Gamma & $[0,\infty)$ & $\frac{\beta^{1-\alpha}}{\Gamma(\alpha)} y^{\alpha-1}\mathrm{e}^{- y/\beta}$ &$1$ & $\frac{N}{2}$ &$ \frac{N^2}{12} $ \\[10pt]
      Pareto & $[\beta,\infty)$ & $\left(\frac{\beta}{y}\right)^{\alpha}$ &$\frac{\alpha-1}{\alpha}$ & $\frac{\alpha N}{2\alpha-1}$ &$ \frac{\alpha^2 (\alpha-1) N^2}{(2\alpha-1)^2(3\alpha-1)} $ \\[15pt]
        \hline
      Distribution & Domain & $h(y)$ as $y \rightarrow y_{\mathrm{max}}$ & $\lambda$ & $E(\tau_N)$ & Var($\tau_N$)\\
      \hline
            &&&&&\\
      Uniform & $[a,b]$ & $\frac{b-y}{b-a}$ &$2$ & $\frac{N}{3}$ &$ \frac{N^2}{18} $ \\[10pt]
      Triangular & $[a,b]$ & $ \frac{(b-y)^2}{(b-a)(b-c)}$ &$\frac{3}{2}$ & $\frac{2N}{5}$ &$ \frac{12 N^2}{175} $ \\[10pt]
            Wigner & $[-R,R]$ & $ \frac{4\sqrt{2}(R-y)^{3/2}}{3\pi R^{3/2}}$ &$\tfrac{5}{3}$ & $\frac{3N}{8}$ &$ \frac{45 N^2}{704} $ \\[10pt]
      Beta & $(0,1)$ & $ \frac{\Gamma(\alpha+\beta)(1-y)^{\beta}}{\Gamma(\alpha)\Gamma(\beta+1)}$ & $\frac{\beta+1}{\beta}$ & $\frac{\beta N}{2\beta+1}$ &$ \frac{\beta^2(\beta+1) N^2}{(2\beta+1)^2(3\beta+1)} $ \\[15pt]
      \hline
    \end{tabular}
  \end{center}
\end{table}

\section{Numerical Comparisons}\label{s.num}

\begin{figure}
\centering
\includegraphics[width=0.475\linewidth]{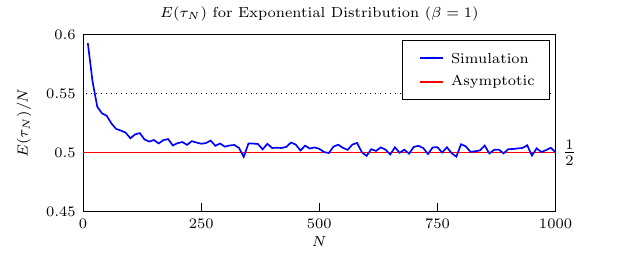}\includegraphics[width=0.475\linewidth]{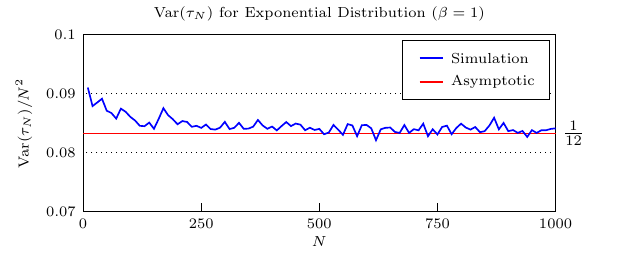}

\includegraphics[width=0.475\linewidth]{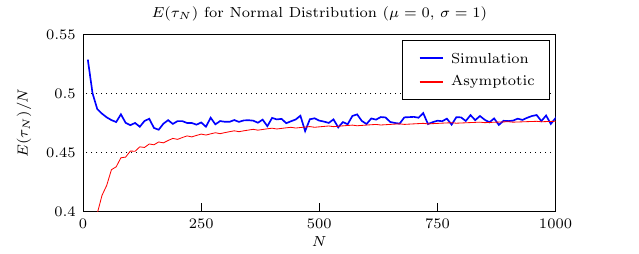}\includegraphics[width=0.475\linewidth]{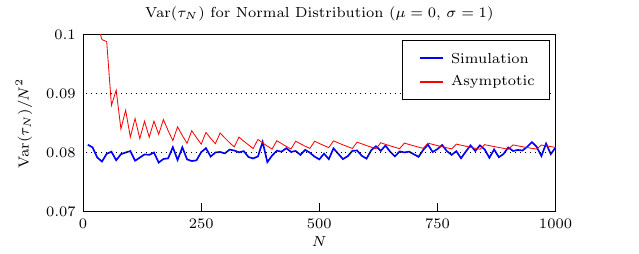}

\includegraphics[width=0.475\linewidth]{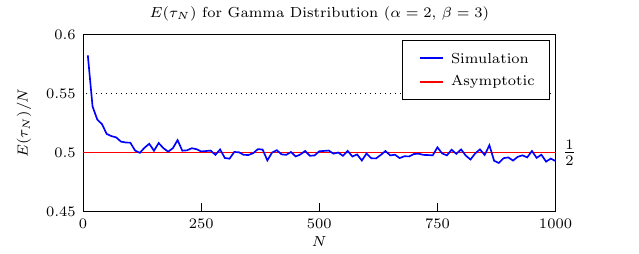}\includegraphics[width=0.475\linewidth]{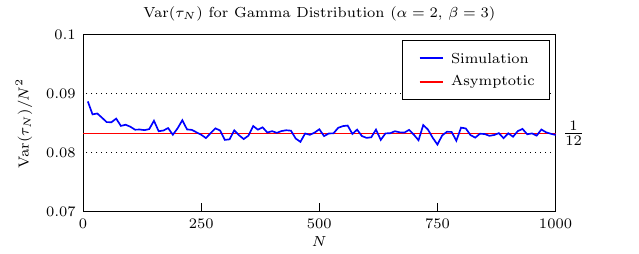}

\includegraphics[width=0.475\linewidth]{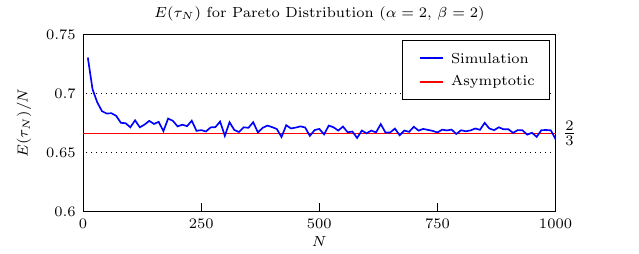}\includegraphics[width=0.475\linewidth]{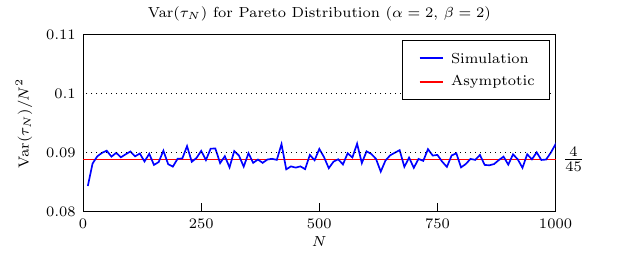}

\includegraphics[width=0.475\linewidth]{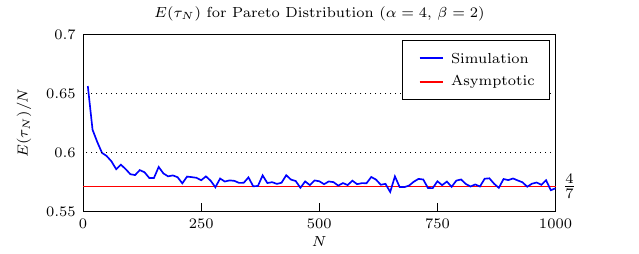}\includegraphics[width=0.475\linewidth]{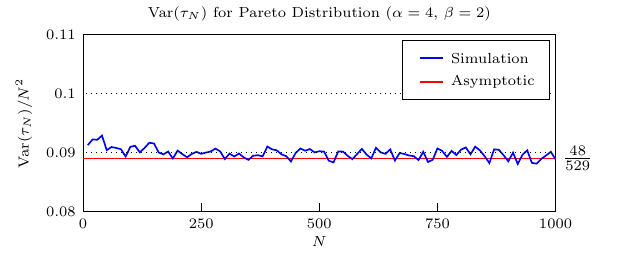}

\includegraphics[width=0.475\linewidth]{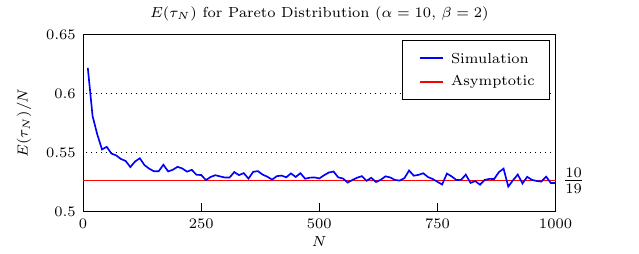}\includegraphics[width=0.475\linewidth]{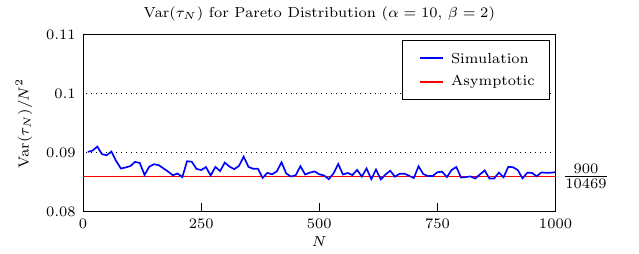}
\caption{Comparison of large-$N$ asymptotic predictions for the expectation and variance of the optimal stopping rule for a range of distributions on unbounded domains. The asymptotic prediction is shown in red, while the simulated result is shown in blue. Each point on the simulated curve was obtained by simulating the optimal stopping problem 10000 times, with the results used to estimate the expectation and variance of the optimal stopping point. We note that, for convenience of presentation, the expectation and variance are scaled by $N$ and $N^2$ respectively, so that the curves tend to a constant value. The first correction term was retained in the asymptotic comparison for the normal distribution, as the asymptotic decay of this term is particularly slow in the large-$N$ limit.}\label{f.unbounded}
\end{figure}

\begin{figure}
\centering
\includegraphics[width=0.475\linewidth]{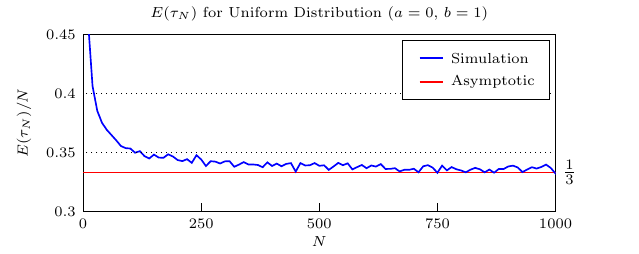}\includegraphics[width=0.475\linewidth]{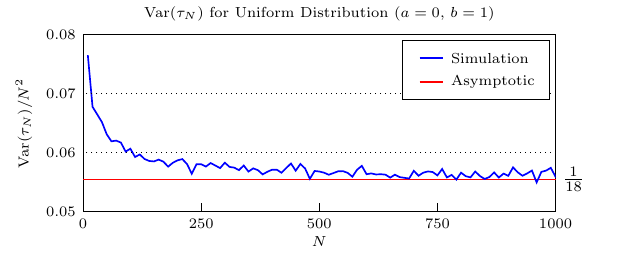}

\includegraphics[width=0.475\linewidth]{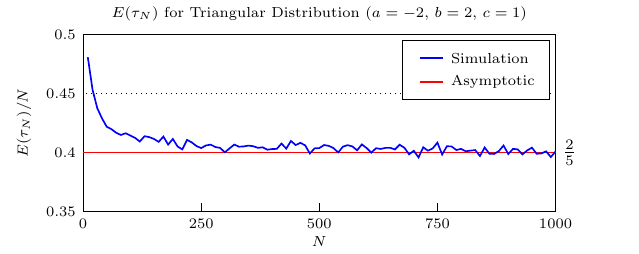}\includegraphics[width=0.475\linewidth]{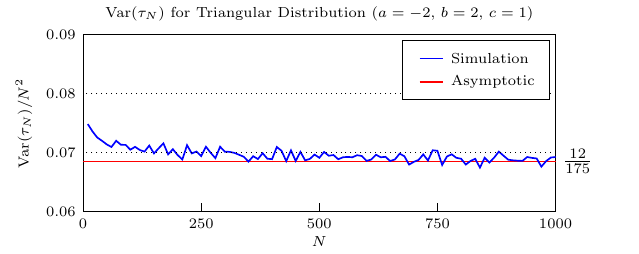}

\includegraphics[width=0.475\linewidth]{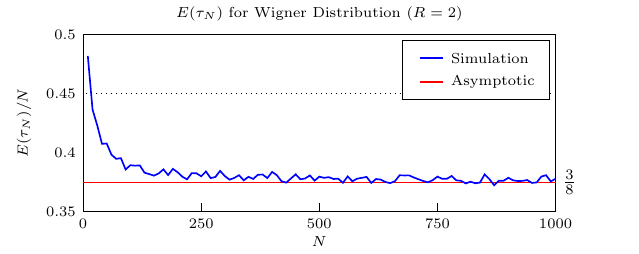}\includegraphics[width=0.475\linewidth]{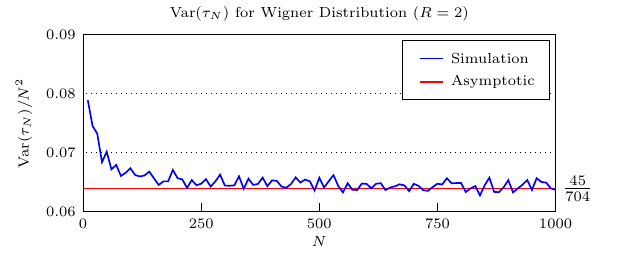}

\includegraphics[width=0.475\linewidth]{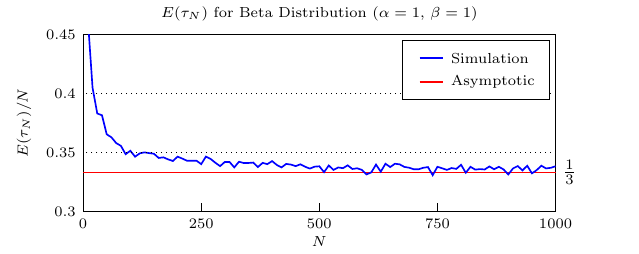}\includegraphics[width=0.475\linewidth]{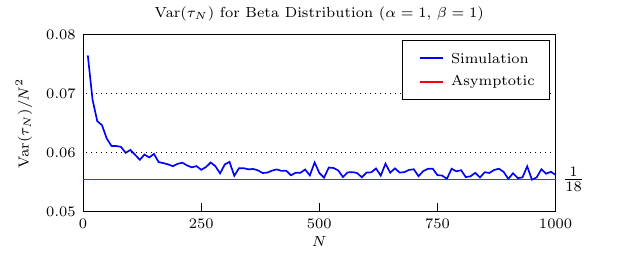}

\includegraphics[width=0.475\linewidth]{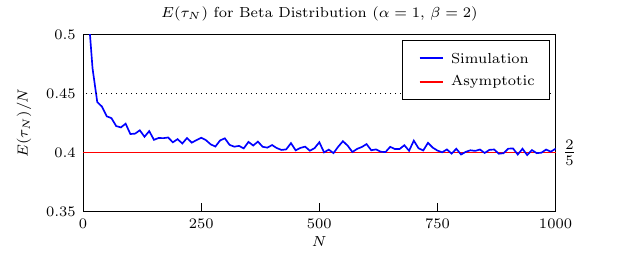}\includegraphics[width=0.475\linewidth]{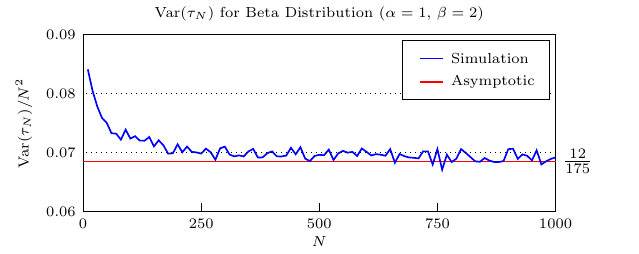}

\includegraphics[width=0.475\linewidth]{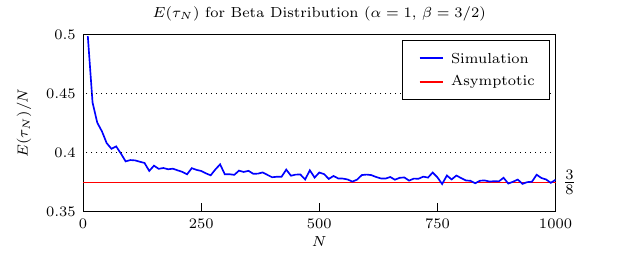}\includegraphics[width=0.475\linewidth]{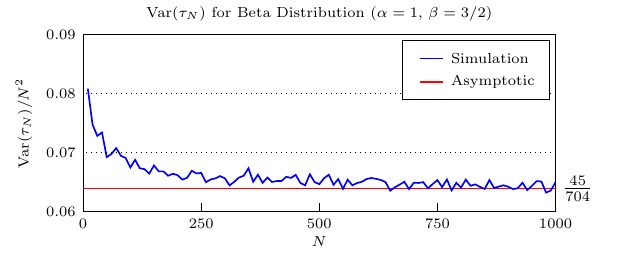}
\caption{Comparison of large-$N$ asymptotic predictions for the expectation and variance of the optimal stopping rule for a range of distributions on bounded domains. The asymptotic prediction is shown in red, while the simulated result is shown in blue. Each point on the simulated curve was obtained by simulating the optimal stopping problem 10000 times, with the results used to estimate the expectation and variance of the optimal stopping point. We note that, for convenience of presentation, the expectation and variance are scaled by $N$ and $N^2$ respectively, so that the curves tend to a constant value.}\label{f.bounded}
\end{figure}

For each of the distributions in Table \ref{T:eqvar}, we validated the asymptotic predictions for the expectation and variance of $\tau_N$ for large $N$ by comparing against numerical simulations. For each value of $N$ from $10$ to $1000$, in increments of ten, we simulated the optimal stopping problem $10000$ times, using the optimal stopping rule from \eqref{1:opstop} to determine the stopping point. The expected value and variance are then estimated for each of these $N$. The results of these computations are illustrated in Figures \ref{f.unbounded} and \ref{f.bounded} for a range of unbounded and bounded distributions respectively.

We note that the calculated expectation and variance are scaled by a factor of $N$ and $N^2$ respectively, so that they tend to constant asymptotic values. The variance of the simulations for these calculated values appears constant under this scaling, which indicates that the simulated expectation and variance both vary more significantly about the asymptotic predicted values as $N$ increases.

In each computation, the numerical results support the corresponding asymptotic predictions. As $N$ increases, the expectation and variance of each set of trials approach the value predicted by the asymptotic calculations. For the Pareto and beta distributions, the numerical comparisons were performed for three different sets of distribution parameters. We see that changing the parameters of the Pareto and beta distributions has the effect of changing the asymptotic expectation and variance, as predicted by the calculated values in Table \ref{T:eqvar}. 

\begin{figure}
\centering
\includegraphics[width=0.475\linewidth]{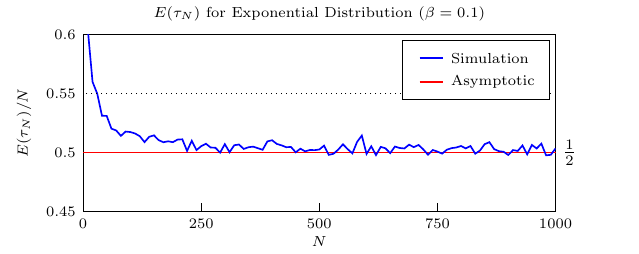}\includegraphics[width=0.475\linewidth]{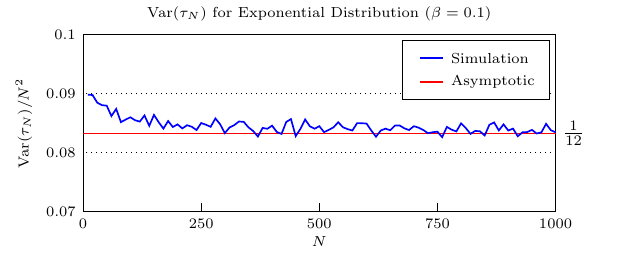}

\includegraphics[width=0.475\linewidth]{E_exp_1.pdf}\includegraphics[width=0.475\linewidth]{var_exp_1.pdf}

\includegraphics[width=0.475\linewidth]{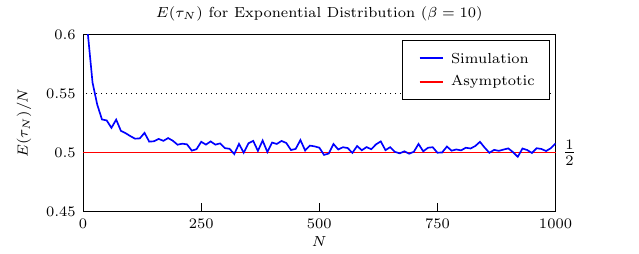}\includegraphics[width=0.475\linewidth]{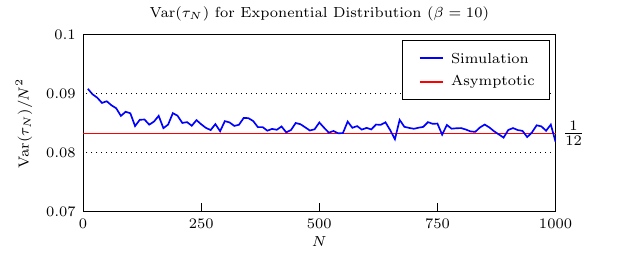}

\caption{Comparison of large-$N$ asymptotic predictions for the expectation and variance of the optimal stopping rule for the exponential distribution with a range of parameter values. In each case, the behaviour tends to the same value in the large-$N$ limit, validating the conclusion that the leading-order asymptotics of statistical properties of optimal stopping do not depend on the distribution parameters for distributions with exponential upper tails. This is in contrast to the Pareto distribution, seen in Figure \ref{f.unbounded}, which has algebraic upper tails; in this case, the asymptotic values change as the parameters are varied.}\label{f.parameter}
\end{figure}

In contrast, Figure \ref{f.parameter} contains simulated results for the exponential function with $\beta = 0.1$, $\beta = 1$, and $\beta = 10$. In each of these figures, the asymptotic behaviour of the expectation and variance of $\tau_N$ tends to identical values for the expectation and variance. These figures support the predictions from Section \ref{s.opstat}, in which it was determined that the leading order large-$N$ approximation of the expectation and variance do not depend on the distribution parameters for distributions with exponential upper tails.

In our asymptotic prediction for the mean and variance of the normal distribution, we retained the first correction term for $h(v_n)$ as $n \rightarrow \infty$, which corresponds to taking the first correction in the expansion of the Lambert-W function \eqref{lambert}. This gives
\begin{equation}
h(v_n) \sim \frac{1}{n} + \frac{1}{n\log(n^2/2\pi)}\quad \mathrm{as} \quad n \rightarrow \infty,
\end{equation}
which allows for a more accurate estimate of $E(\tau_N)$ using \eqref{3:Nsum} and $\var(\tau_N)$ using \eqref{3:Nsum0}. This first correction is small compared to the leading order behaviour as $N \rightarrow \infty$, but decays so slowly that it must be included in order to obtain accurate predictions for the expectation and variance of $\tau_N$ for even moderately large values of $N$. By comparing the behaviour of the simulated and asymptotic results, we see that the computations agree with this more accurate asymptotic prediction.

%
%
%
%
%
%
%
%
%

\section{Conclusions}

In this paper we have derived asymptotics of optimal stopping times for sequences of independent identically distributed continuous random variables. In particular, we have found the asymptotics of the expected value and the variance of the stopping time for large classes of density functions whose domains have either infinite or finite upper bounds. 


Asymptotic calculations were performed on a number of probability distributions, on both bounded and unbounded domains. Numerical simulations were subsequently performed to calculate the expectation and variance of the optimal stopping rule for a range of values of $N$, ranging from 10 to 1000. In each case, the simulated results tended towards the asymptotic prediction in the large-$N$ limit, validating the asymptotic approach. 


One particularly interesting observation is that, if a density function $f(y)$ has no upper bound and the upper tail decays exponentially as $y \rightarrow \infty$, the expectation and variance of the optimal stopping rule are given by $E(\tau_N) = N/2$ and $\mathrm{Var}(\tau_N) = N^2/12$. These asymptotic values do not depend on any other features of the distribution, and hold for any distribution with exponentially-decaying upper tail behaviour. This parameter independence was tested for the exponential distribution, for which the simulated values converged to the large-$N$ asymptotic prediction given a range of different parameter choices.

The independence of the expectation and variance of the optimal stopping rule from distribution parameters does not hold for distributions with algebraically-decaying upper tails, such as the Pareto distribution. We found in our analysis that the expectation and variance of the optimal stopping rule in this case depends on the algebraic power of the tail decay rate, which was supported by numerical computation for a range of different parameter values.

In most cases, the leading-order asymptotic behaviour was sufficient to explain the simulated results; however, this is not necessarily true. Comparing the leading-order asymptotic approximation with the simulated results for the normal distribution showed significant error between the predicted and simulated values. This error was resolved by including an extra correction term in the asymptotic calculations, leading to agreement between the asymptotic and simulated results. It is important to note that for some distributions, higher correction terms to the value of a sequence can be required in order for the asymptotics to be useful in predicting statistical properties for moderately large values of the asymptotic parameter.





The key contribution from this study is outlining a general strategy for obtaining the asymptotic behaviour of the expectation and variance of the optimal stopping rule that can be applied to a wide range of distributions in a consistent fashion. This method recovers existing results from \cite{mazalov2004asymptotic}, while also making predictions for many other distributions that were subsequently validated through comparison with simulated results.


\appendix

\bibliographystyle{plain}
\bibliography{StoppingTimes}

\end{document}